\documentclass[11pt]{article}
\usepackage{graphicx}
\usepackage{endnotes}
\usepackage{hyperref}
\usepackage{natbib}
\usepackage{amsmath,amsfonts,amsthm,bm}

\let\footnote=\endnote

\def\argmin{\mathop{\rm arg\,min}\limits}

\def\R{\ifmmode{I\hskip -3pt R}
           \else{\hbox{$I\hskip -3pt R$}}\fi}

\def\N{\ifmmode{I\hskip -3pt N}
           \else{\hbox{$I\hskip -3pt N$}}\fi}
\def\m{{\bf m}}
\def\x{{\bf x}}

\def\z{{\bf z}}
\def\X{{\bf X}}
\def\Y{{\bf Y}}
\def\Z{{\bf Z}}
\def\z{{\bf z}}
\def\G{{\bf G}}
\def\I{{\mathbf I}}
\def\bmu{{\bm{\mu}}}
\def\bS{{\bm{\Sigma}}}

\newcommand{\Cor}{{\mathrm{Cor}}}

\newtheorem{theorem}{Theorem}[section]
\newtheorem{definition}[theorem]{Definition}
\newtheorem{lemma}[theorem]{Lemma}

\newtheorem{example}[theorem]{Example}
\newtheorem{remark}[theorem]{Remark}
\newtheorem{corollary}[theorem]{Corollary}
\title{{\bf Parameters not empirically identifiable or distinguishable, including correlation between Gaussian observations}}
\author{Christian Hennig,
 Dipartimento di Scienze Statistiche ``Paolo Fortunati'',\\ Universita di Bologna,\\
 Via delle Belle Arti, 41, 40126 Bologna\\
 christian.hennig@unibo.it}

\begin{document}
\maketitle

\noindent
{\bf Abstract:}~\\
{\bf Note} Accepted version, published in {\it Statistical Papers},\\ \verb|https://doi.org/10.1007/s00362-023-01414-3|.
~\\~\\
It is shown that some theoretically identifiable parameters cannot
be empirically identified, meaning that no consistent estimator of them can exist.
An important example is a constant correlation between Gaussian observations (in presence of
such correlation not even the mean can be empirically identified). 
Empirical identifiability and three versions of empirical distinguishability 
are defined. 
Two different constant
correlations between Gaussian observations cannot even be empirically distinguished. A further example are cluster membership parameters in $k$-means clustering. Several existing results in the literature are connected to the new framework. General conditions are discussed under which independence can be distinguished from dependence.

\section{Introduction}\label{sintro}
Meaningful statistical inference is only possible if the target of inference 
(parameter) is identifiable, meaning that if parameter values differ,
the parameterised distributions should also differ. There are several 
versions of identifiability definitions, and many identifiability and non-identifiability results, see, e.g., \cite{YakSpra68,Rothenberg71,PrakasaRao92,HoRos17}. 

Here situations are treated in which parameters are identifiable according to 
this classical definition, yet the parameters cannot be identified
from observed data. 
\cite{Rothenberg71,Hsiao83,PrakasaRao92} define identifiability
with explicit reference to observable data but do not cover the 
issues that are treated here. Regarding the results in
\cite{Rothenberg71}, there is no difference between classical identifiability 
and identifiability from observations. Simple examples for classical 
identifiability issues are the non-identifiability of linear regression
parameters in case of collinear explanatory variables, and identifiability
of the parameters of mixture distributions, which can be guaranteed under 
certain assumptions (\cite{YakSpra68}), particularly ruling our label switching,
but counterexamples exist \cite[Chapter 8]{PrakasaRao92}.  
\cite{Hsiao83,PrakasaRao92} also study
situations in which issues occur because certain modelled 
random variables are unobservable, 
such as the true value of a variable in errors-in-variables 
models. Some examples in Section \ref{sliterature} are also of this kind, but
there are further reasons why the observed data may not allow for 
identification of classically identifiable parameters, which are explored 
here.

Some such situations have already appeared in the literature, 
see, e.g., \cite{NeySco48,BahSav56,Donoho88,SpGlSc93,RSSW03,MBSK08,AlmMou14}.
Section \ref{sliterature} gives more details on these works, 
and how they fit into the unified terminology introduced here.

It turns out that there are different possible levels of information about identifiable parameters in the data, and therefore various definitions are introduced. Consistent estimators may or may not exist
(``empirical identifiability''). Sets that can distinguish two parameter values 
for a finite sample size may or may not exist 
(``empirical distinguishability'' with weaker and stronger versions).  

The concept of empirical identifiability is closely connected to the concept of estimability, which is also stronger than classical identifiability. Once more there are several versions around. The concept mostly focuses on what can be estimated with a give finite sample, and its connection with classical identifiability is investigated, see, e.g., \cite{BunBun74,JacGre85,MacNic20}. 

This work was motivated by the discovery that data hold no information about distinguishing i.i.d. Gaussian observations from Gaussian data with a constant correlation between any two observations. This will be used as a guiding example. 
Section \ref{sgausscor} derives a key result regarding this situation. Section \ref{sidentdist} presents the main definitions, some of their implications and some more examples. Section \ref{sliterature} reviews results from the literature that fit into the framework of Section \ref{sidentdist}. Section \ref{sindep} uses this framework to discuss the general problem of telling apart dependence and independence in situations in which potential dependence is not governed by the observation order or observable external information. This is relevant in many situations that require independence assumptions. Section \ref{skmeans} presents another example in some detail, namely the empirical identification of parameters indicating the cluster memberships of every single point in $k$-means clustering. Section \ref{sconc} concludes the paper. All proofs are in the Appendix.

\section{Constant correlation between Gaussian observations} 
\label{sgausscor}
This work was motivated by the following example, which in itself should be 
of strong interest.
\begin{example}\label{exgausscor}
A model assumption for much standard statistical inference is to assume
independently identically distributed (i.i.d.) 
Gaussian $X_1,\ldots,X_n,\ X_1\sim{\cal N}(\mu,\sigma^2)$ (model M0). Now
consider Gaussian $X_1,\ldots,X_n$ with 
correlation $\Cor(X_i,X_j)=\rho>0$ constant for any $i\neq j$ (model M1; M01 denotes the model with $\rho\ge 0$ assumed).

This would be a problem for inference about $\mu$, because in the latter
situation, for the arithmetic mean $\bar{X}_n$:
\[
{\cal L}(\bar{X}_n)={\cal N}\left(\mu,\frac{(1-\rho)\sigma^2}{n}+\rho\sigma^2\right)\to_{n\to\infty} {\cal N}\left(\mu,\rho\sigma^2\right).
\]
This means that the mean is inconsistent for $\mu$ as long as $\rho\sigma^2>0$;
confidence intervals and tests computed based on the i.i.d. assumption will be
biased, possibly dramatically so.  

Although $\rho$ is identifiable in the classical sense, 
it turns out to not be empirically identifiable from observed data.
It is not even possible to empirically distinguish any two $\rho_1\neq\rho_2$ (in fact, plots of data generated from model M01 with different values for $\rho$ including the i.i.d. case $\rho=0$ do not reveal any features by which these distributions could be distinguished). 
$\mu$ in model M01 is not empirically identifiable either, but two $\mu_1\neq\mu_2$ are empirically distinguishable, see Section \ref{sidentdist}.
\end{example}

The following lemma shows that in model M01, the conditional
distribution given the mean $\bar{X}_n$ is the same as for i.i.d., therefore 
uncorrelated, Gaussian random variables. But the mean does not hold information 
about correlations (or rather, any information about correlations is confounded
with the information about the true means), meaning that model M1 cannot be 
distinguished from model M0 based on the data alone.
\begin{lemma}\label{lcondgauss} 
$\X_n=\begin{pmatrix} X_1\\ \vdots \\ X_n\end{pmatrix},\ \Y_n=\begin{pmatrix} Y_1\\ \vdots \\ Y_n\end{pmatrix}$. Assume 
\[ {\cal L}(\X_n)
= {\cal N}_n(\bm{\mu,\Sigma}),\ \bmu=
\begin{pmatrix}\mu\\ \vdots\\ \mu\end{pmatrix},\
\bS= 
\begin{bmatrix} 
    \sigma^2 & \rho\sigma^2 & \dots & \rho\sigma^2\\
    \rho\sigma^2 & \sigma^2 & \dots & \rho\sigma^2\\
    \vdots & \vdots & \ddots & \vdots\\
    \rho\sigma^2 & \rho\sigma^2 &\dots  & \sigma^2 
    \end{bmatrix}.
\]

Then, for 
\[
{\cal L}(\Y_n)={\cal N}_n(\bmu,(1-\rho)\sigma^2\I_n):\
{\cal L}(\X_n\mid \bar{X}_n)={\cal L}(\Y_n\mid \bar{Y}_n),
\]
which does not depend on $\mu$.
\end{lemma}
Thus, conditionally on the mean, $\X_n$ will look like i.i.d. Gaussians with
variance $(1-\rho)\sigma^2$; for $\rho>0$ there is less variation of
the $X_i$ given their mean than their unconditional variance. On the other
hand, $\bar{X}_n$ has a larger variance than under independence ($\rho=0$).

In fact, the model can equivalently be written as a model with a 
single realisation of a random effect $Z,\ i=1,\ldots,n$:
\begin{equation}\label{erandomeff}
X_i=\mu+Z+E_i,\ Z\sim{\cal N}(0,\tau_1^2),\ E_i\sim{\cal N}(0,\tau_2^2),\ 
\sigma^2=\tau_1^2+\tau_2^2,\ \rho=\frac{\tau_1^2}{\tau_1^2+\tau_2^2},
\end{equation}
which suggests that observed data look like i.i.d. ${\cal N}(\mu^*,\tau_2^2)$
with $\mu^*=\mu+Z$, and $Z$ is unknown and unobservable. The definitions and results in Section \ref{sidentdist} aim at making precise a general sense in which observations give no information about $\rho$ and limited information about $\mu$.

\section{Empirical identifiability and distinguishability}\label{sidentdist}
Let $X_1,X_2,\ldots,X_n,\ldots$ be random variables on a space ${\cal X}$,
for $n\in\N:\ {\cal L}(X_1,\ldots,X_n)=P_{n;\theta}$ with parameter 
$\theta\in \Theta$. $P_{\infty;\theta}$ denotes the distribution of the whole sequence.
The spaces ${\cal X}$ and $\Theta$ can be very general, but
assume that $\Theta$ is a metric space with metric $d_\Theta$.
The focus may be on the parameter $\theta$ in full, or it may be on
$g(\theta)$, where 
$g:\ \Theta\mapsto\Lambda$, $\Lambda$ being a metric space with metric 
$d_\Lambda$. No further conditions on $\Theta$ and $\Lambda$ are
required for the general definitions. It is generally assumed that the 
underlying $\sigma$-algebras are rich enough so that the sets required in
the arguments are measurable. In the specific cases discussed here
this is always fulfilled using standard (Borel) 
$\sigma$-algebras and parameter spaces.
Sometimes but not always $g(\theta)=\theta$
and $\Lambda=\Theta$ are considered. Other examples for $g$ are a 
projection on a lower dimensional space, or  
an indicator function for a parameter subset (hypothesis) of interest. 
\begin{definition}
$g(\theta)$ is called {\bf empirically identifiable} if it is possible to 
find a consistent sequence of estimators
$(T_n)_{n\in\N}$, i.e., with $T_n={\cal X}^n\mapsto \Lambda$,
$\forall \theta\in\Theta:\ T_n(X_1,\ldots,X_n)\to g(\theta)$ in probability.
\end{definition}
Consistency is always meant with respect to $P_{\infty;\theta}$.
Traditionally, statistical identifiability of a parametric model 
$\left(P_\theta\right)_{\theta\in\Theta}$ 
means that $\theta_1\neq \theta_2\Rightarrow P_{\theta_1}\neq P_{\theta_2}$;
for parameter parts, 
$g(\theta_1)\neq g(\theta_2)\Rightarrow P_{\theta_1}\neq P_{\theta_2}$ is often referred to as partial identifiability (\cite{PrakasaRao92,HoRos17}). If parameters are not (partially) identifiable, they can obviously not be empirically identifiable, because
no consistent estimator can tell equal distributions apart: 
\begin{corollary} Parameters and parameter parts that are empirically identifiable
are also identifiable.
\end{corollary} 
Here, data generating mechanisms are treated that do not allow to empirically 
identify 
parameters that are in fact identifiable in the traditional sense.
Model M01 is an example. Obviously, distributions
with different correlation parameters $\rho_1\neq\rho_2$ are different from 
each other, and $\rho$ can be estimated consistently if the whole sequence of 
$n$ observations is repeated {\it independently}. 
In this case, assuming equal correlation between any two components, 
the sequence of length $n$ of observations becomes an $n$-variate 
Gaussian, and the sample correlation between any two of the $n$ components 
will estimate $\rho$ consistently,  
although
a better estimator will of course use information from all components.
The data generating mechanism
modelled in Section \ref{sgausscor} does not allow for independent repetition;
all available observations are dependent from all other observations, and this
makes consistent estimation of $\rho$ impossible:
\begin{theorem}\label{tnotidrho} Using the notation of Lemma \ref{lcondgauss}, if, 
for $n\in\N:\ {\cal L}(\X_n)= {\cal N}_n(\bm{\mu,\Sigma})$, then, with 
$\theta=(\mu,\sigma^2,\rho)$, 
$g(\theta)=\rho$ is not empirically identifiable in model M01. 
\end{theorem} 

Not only is the correlation $\rho$ not empirically identifiable, the same holds 
for $\mu$, meaning that in practice using an estimator different from 
$\bar{X}_n$ does not help dealing with the potential existence of $\rho> 0$.

\begin{theorem}\label{tmunonid} 
Using the notation of Lemma \ref{lcondgauss} and Theorem \ref{tnotidrho}, if 
for $n\in\N:\ {\cal L}(\X_n)= {\cal N}_n(\bm{\mu,\Sigma})$, then 
$g(\theta)=\mu$ is not empirically identifiable. 
\end{theorem} 
The proof of Theorem \ref{tmunonid} relies on the random effects formulation
(\ref{erandomeff}) with only a single realisation of the random effect. 
A similar case can be made for a standard random effects 
model assuming that the number of realised values of the random effect is 
bounded even if the number of observations goes to infinity, 
as expressed in the following model M2:
\begin{displaymath} 
X_{ij}=\mu+Z_i+E_j,\ Z_i\sim{\cal N}(0,\tau_1^2),\ E_j\sim{\cal N}(0,\tau_2^2),
\end{displaymath}
$i=1,\ldots,m$ (group), $j=1,\ldots,n_i$ (within group observation),
$n=\sum_{i=1}^m n_i,\ \theta=(\mu,\tau_1,\tau_2)$. Let $m$ be fixed, whereas $n$ is allowed to grow. Let $\X_n$ be the vector collecting all $X_{ij}$.

Such a model 
could make sense for a random effects meta analysis with a low number $m$
of studies, each of which is potentially large, but it does not allow to
empirically identify the random effects' variance, and neither the overall 
mean, unless $m\to\infty$. Consequently, common advice in
the meta analysis literature is to not use a random effects model if the number
of studies is low (see, e.g., \cite{KuMoSt08}). 
\begin{lemma}\label{lrandomeff}
In model M2, $g_1(\theta)=\mu$ and $g_2(\theta)=\tau_1^2$ are not
empirically identifiable, whereas 
$g_3(\theta)=\tau_2^2$ is empirically identifiable.
\end{lemma}
In model M01, there is a difference between trying to estimate $\rho$ on one hand
and $\mu$ on the other hand. While $\mu$ cannot be estimated consistently, 
in case that $\rho$ is small, the data can give fairly precise information about
its location, whereas there is no information in the data about $\rho$ at all.
The following definition aims at formalising this difference.
\begin{definition} \label{dindistset} For $n\in\N,\  \alpha\le\beta\in(0,1]$,
an observable set $A$, i.e., any measurable set expressing an 
observable event,
is an {\bf $(\alpha,\beta,n)$-distinguishing set} for 
$\theta_1\neq\theta_2\in\Theta$ if 
\begin{equation}\label{eindistset}
P_{n;\theta_1}(A)\le\alpha,\ P_{n;\theta_2}(A)>\beta.
\end{equation}
\end{definition} 
\begin{definition} \label{dindist}
Two values $\lambda_1\neq\lambda_2\in\Lambda$ 
are called 
{\bf empirically distinguishable} if $\exists n, \alpha\in(0,1],$ and $\forall \theta_1,\theta_2\in\Theta \mbox{ with }g(\theta_1)=\lambda_1,\ 
g(\theta_2)=\lambda_2$ there is an $(\alpha,\alpha,n)$-distinguishing set $A$. 
\end{definition} 
Obviously, this definition is symmetric in $\lambda_1, \lambda_2$. Before
returning to the problem of constant correlation between Gaussian observations,
empirical distinguishability is discussed in some more generality.  

For $\epsilon>0$ and $\eta_0$ in some metric space $H$ with metric $d_H$, define $B_\epsilon(\eta_0)=\{\eta:\ d_H(\eta,\eta_0)\le\epsilon\}$.
If $g(\theta)=\theta$, empirical distinguishability follows from empirical identifiability, because there is a consistent estimator $T_n$ of $\theta$, and 
$A=\{T_n\in B_\epsilon(\theta_2)\}$ will distinguish $\lambda_1=\theta_1\neq\lambda_2=\theta_2$ for large enough $n$ if $\epsilon$ is chosen small enough that 
$\theta_1\not\in B_\epsilon(\theta_2)$. 

In general, empirical identifiability does not imply empirical distinguishability. If $g(\theta)$ specifies only a part of the information in $\theta$, it may happen that no set $A$ can distinguish $g(\theta_1)=\lambda_1$ from $g(\theta_2)=\lambda_2$ uniformly over the information in $\theta$ that is not in $g(\theta)$, even if $g(\theta)$ is empirically identifiable.

\begin{example}\label{exindist} 
Let $X_i,\ i\in\N$ be independently distributed according 
to $P_\theta,\ \theta=(p,m),\ p\in [0,1],\ m\in\N,$ defined as follows: 
For $i\ge m$, 
${\cal L}(X_i)=$Bernoulli$(p)$. For $i<m$, ${\cal L}(X_i)=$Bernoulli$(q)$, where
$q$ is randomly drawn from ${\cal U}(0,1)$. $g(\theta)=p$ is empirically identifiable, because $\bar X_n$ is consistent for it. But any $p_1\neq p_2$ are not empirically distinguishable, because for any $n<m$, $X_1,\ldots,X_n$ do not contain any information about $p$. In this situation, the data may carry information about whether $n>m$ (namely where it can be observed that at some point in the past $q$ may likely have changed to $p$), in which case it also carries information about $p$, but for $m=1$ this can never happen, and for very small $m$ this can hardly ever be diagnosed with any reliability. 
\end{example}

If $\lambda=g(\theta)$, in order to make $\lambda_1$ and 
$\lambda_2$ empirically distinguishable from having a consistent estimator 
$(T_n)_{n\in\N}$
of $\theta$, in general $g$ needs to be uniformly 
continuous, and $(T_n)_{n\in\N}$ needs to be
uniformly consistent on $C=g^{-1}(\lambda_1)\cup g^{-1}(\lambda_2)$, i.e.,
 \[
\forall \epsilon>0,\alpha>0\ \exists n_0 \forall n\ge n_0, 
\theta\in C:\
P_{n;\theta}\{T_n\in B_\epsilon(\theta)\}>1-\alpha.
\]   
The latter holds automatically in case that $\theta$ is empirically 
identifiable if $|C|$ is finite, and in particular if $g$ is bijective.  
\begin{lemma} \label{lindist}
If $\theta\in\Theta$ is empirically identifiable,
$g$ is uniformly 
continuous on an open superset of $C$,
requiring that such a set exists, and there is an 
estimator $T_n$ of $\theta$ that is uniformly consistent on  
$C$, then
$\lambda_1\neq\lambda_2\in\Lambda$ are empirically distinguishable. 
\end{lemma}
Only assuming $\lambda=g(\theta)$ but not $\theta$ to be empirically identifiable, consistency of $(T_n)_{n\in\N}$ as estimator of $g(\theta)$ needs to be uniform on $C$ for making $\lambda_1\neq\lambda_2$ empirically distinguishable.

Consistent estimators are uniformly consistent in many situations. For example,
the behaviour of affine equivariant multiple linear regression estimators 
is uniform over the whole parameter space. Therefore, if they are consistent for the full parameter vector, they are uniformly consistent, and will 
according to Lemma \ref{lindist} empirically 
distinguish subvectors and single coefficients of the regression parameter.

\begin{remark}\label{rindist} There are possible variants of Definition \ref{dindist} that may be taken as different ``grades'' of empirical distinguishability.

In the situation in Example \ref{exindist}, events may be observed 
for large enough $n$ that can distinguish $p_1$ and $p_2$, even though this is 
not guaranteed to happen. A concept that allows $p_1$ and $p_2$ to be seen
as distinguishable in some sense (at least if too low $m$ is excluded; $m$ needs to be large enough that a ``change point'' after $m$ observations can be 
diagnosed with nonzero probability regardless of $p$)
is ``potential distinguishability'', see
Definition \ref{dmoredist}. 

Furthermore, it would make a difference to not allow 
the probabilities  $P_{n;\theta_1}(A),\ P_{n;\theta_2}(A)$ to be arbitrarily close.
Consider a simple i.i.d. ${\cal N}(\theta,1)$ model. With the given definition,
and, for given fixed $\theta_0$, $g(\theta)=1(\theta=\theta_0)$, 0 and 1 can be empirically distinguished (the rejection region of the standard two-sided test will do), i.e., it can be distinguished whether $\theta=\theta_0$ or
$\theta\neq \theta_0$. Using a definition that requires a ``gap'' of some 
$\beta>0$ between
$P_{n;\theta_1}(A)$ and $P_{n;\theta_2}(A)$,
0 cannot be distinguished from 1, because for $\theta\to\theta_0$, $P_\theta(A)$ for any $A$ gets arbitrarily close to $P_{\theta_0}(A)$. Both of these definitions can be seen as appropriate, from different points of view. 
It could be argued that $\bar X$ contains some, if not necessarily conclusive, information about whether $\theta=\theta_0$, and it should therefore count as distinguishable from $\theta\neq\theta_0$, as achieved by 
Definition \ref{dindist},

But even with arbitrarily large $n$, $\bar X$ will not be exactly zero, and will therefore be at least as compatible with some $\theta\neq\theta_0$ as with $\theta_0$, which could be used to argue that the two should not be defined as distinguishable. Choosing $g$ as the identity, any fixed $\theta\neq\theta_0$ could still be distinguished from $\theta_0$; in this case there is a positive distance between $\theta$ and $\theta_0$, which makes $\bar X_n$ a better fit for the closer parameter. This can be achieved by the concept of ``empirical gap distinguishability'', see Definition \ref{dmoredist}.

\end{remark}

\begin{definition} \label{dmoredist}
\begin{enumerate}
\item[(a)]
$\lambda_1\in\Lambda$ is called 
{\bf potentially empirically distinguishable} from $\lambda_2\in\Lambda$
if $\exists \alpha\in(0,1),\ n\in \N$, a set $D\subseteq g^{-1}(\lambda_1)\times g^{-1}(\lambda_2)$ so that
\[
\forall \theta_1\in g^{-1}(\lambda_1)\ \exists \theta_2\in g^{-1}(\lambda_2):\
(\theta_1,\theta_2)\in D,
\]
and a set $A$ that $(\alpha,\alpha,n)$-distinguishes $\theta_1$ from $\theta_2$ for all $(\theta_1,\theta_2)\in D$.
\item[(b)] Two values $\lambda_1\neq\lambda_2\in\Lambda$ 
are called 
{\bf empirically gap distinguishable} if $\exists n, \alpha<\beta\in(0,1],$ and $A$
that is an $(\alpha,\beta,n)$-distinguishing set  
$\forall \theta_1,\theta_2\in\Theta \mbox{ with }g(\theta_1)=\lambda_1,\ 
g(\theta_2)=\lambda_2$.
\end{enumerate}
\end{definition} 
Potential distinguishability in particular implies $\forall \theta_1\in g^{-1}(\lambda_1):\ P_{n;\theta_1}(A)\le\alpha$, so that $\lambda_1$ can be 
``rejected'' by the indicator of $A$ (keeping in mind that $\alpha$ cannot 
necessarily be chosen small), even though this test may
be biased against some $\theta_2\in g^{-1}(\lambda_2)$.
Potential distinguishability is not symmetric in $\lambda_1$ and $\lambda_2$, but in Example \ref{exindist}, in fact $p_1$ is potentially distinguishable from
$p_2$ ($A$ can be chosen as intersection of a set that rejects the null 
hypothesis of no change point in the binary sequence, see \cite{Worsley83}, 
and $\lvert\bar X_n^*-p_1\rvert$ being large, where $\bar X_n^*$ is the mean after the estimated change point), and $p_2$ is potentially 
distinguishable from $p_1$ in the same way. The reason why there is symmetric
potential distinguishability here but not standard distinguishability is that
different distinguishing sets are required for the two directions, and 
that no set works uniformly over all $\theta\in g^{-1}(p_1)\cup g^{-1}(p_2)$. See Example \ref{edensity} for a genuinely asymmetric instance of potential distinguishability.

Empirical gap distinguishability is stronger than empirical distinguishability,
whereas potential distinguishability is weaker.  
Still, lack of empirical gap distinguishability means that in terms of the
parameterised probabilities of observable sets, 
$\lambda_1\neq\lambda_2$ appear 
infinitesimally close, regardless of $d_\Lambda(\lambda_1,\lambda_2)$.

Theorem \ref{tgaurhoindist} treats empirical distinguishability in model M01.
\begin{theorem}\label{tgaurhoindist}
In model M01, any two $\rho_1\neq\rho_2\ge 0$,
$\rho_1,\rho_2<1$, are indistinguishable and not even
potentially distinguishable in any direction,
whereas any two $\mu_1\neq \mu_2$ are empirically distinguishable.
\end{theorem}

\begin{remark}
Apart from partial identifiability, there are further weaker 
versions of the classical identifiability concept. A parameter value is locally 
identifiable if in an open neighbourhood in the parameter space there is no other
parameter that parameterises the same distribution (\cite{Rothenberg71}). Set
identifiability (\cite{HoRos17}) 
means that sets of equivalent parameter values can be identified and potentially
be estimated
as opposed to a precise parameter value. In other words, parameters from non-equivalent sets could be (empirically) 
distinguished, whereas equivalent parameters could not be distinguished. 
Empirical versions of these definitions
are possible, but all the situations treated here that are not empirically identifiable would not be empirically locally or set identifiable either. All proofs of empirical non-identifiability in the Appendix rule out the existence of consistent estimators that can tell any two parameter values apart, so neither 
parameter sets nor parameter values in any neighbourhood of each other
can be consistently told apart.
\end{remark}

\section{Examples from the literature}\label{sliterature}
The concepts of empirical identifiability and distinguishability provide a framework that fits various existing results on the limitations of empirically identifying parameters or hypotheses that are identifiable according to the classical definition. 
\begin{example}\label{eneysco} The so-called ``incidental parameter problem'' was introduced by \cite{NeySco48}, see \cite{Lancaster00} for a review. It refers to a situation in which there are observed units the number of which is allowed to go to infinity, and for these observed units there is a bounded finite number of observations, say $x_{ij},\ i=1,\ldots,n,\ j=1,\ldots,t$, where $n\to\infty$ but $t$ fixed. The model for the distribution of the corresponding random variables $X_{ij}$ involves some parameters $\alpha_i,\ i=1,\ldots,n$. For the estimation of each of these there are only $t$ observations available, and the $\alpha_i$ will not be estimated consistently if $n\to\infty$, so they are not empirically identifiable, although they may well be empirically distinguishable, depending on the specific model. The problem can be avoided in many situations by modeling the $\alpha_i$ as random effects, so that they are governed by only one or few parameters, but in some situations, e.g., panel data in economics (\cite{Lancaster00}), researchers may be interested in inference about specific $\alpha_i$, and also standard distributional assumptions for the random effects distribution may not seem realistic. 

\end{example}

\begin{example}\label{ebahsav} In a famous paper, \cite{BahSav56} show the non-existence of 
valid statistical inference for the problem of finding out about the true mean in a sufficiently large class of distributions ${\cal F}$ with existing mean (essentially requiring that for every $P$ and $Q$ also any mixture of them is in ${\cal F}$, and that
$\forall \mu\in\R\ \exists P\in{\cal F}:\ E_P(X)=\mu$). 

Applying the terminology of the present paper, their Theorem 1 shows that any two means $\mu_1$ and $\mu_2$ are not empirically gap distinguishable. The reason is that with $E_P(X)=\mu_1$, $E_Q(X)$ can be chosen so that for arbitrarily small $\epsilon>0$ and $R=(1-\epsilon)P+\epsilon Q$, $E_R(X)$ can take any value $\mu_2\neq\mu_1$. For arbitrarily large $n$ and small enough $\epsilon$, $P^n(A)$ and $R^n(A)$ are arbitrarily close.
\end{example}

\begin{example}\label{edensity} Expanding the work of \cite{BahSav56}, \cite{Donoho88} 
considered functionals $J$ of distributions, including the number of modes
of the density, the Fisher information, any $L_p$-norm of any derivative of the 
density, the number of mixture components, and the negentropy. He showed 
that for a sufficiently rich nonparametric family ${\cal P}$ of distributions,
the graph (i.e., the set of pairs $(P,J(P))$ with distribution $P$ and 
functional value $J(P)$) is dense
in the epigraph (the set of pairs $(P,j)$ where $j\ge J(P)$) using a
``testing topology'' induced by $d(P,Q)=\sup_{0\le\psi\le 1} 
\lvert\int\psi dP-\int \psi dQ\rvert$. As $\psi$ can be the indicator of a 
distinguishing set, two classes of distributions ${\cal P}$ and ${\cal Q}$ 
are not empirically gap distinguishable if $\inf_{P\in{\cal P}, Q\in{\cal Q}}d(P,Q)=0$.
On the positive side, Donoho proves the existence of confidence sets for lower
bounds of these functionals (except the negentropy); enough data can make it 
possible to identify that $J(P)\ge k$ fixed and given, provided that this is
indeed the case for $P$. $J(P)=k_1$ can be potentially distinguished from
$J(P)=k_2<k_1$, but not from $J(P)=k_3>k_1$.

Whereas distinguishability results are mostly negative, \cite{Donoho88} 
constructs a consistent estimator of the number of modes (Corollary of 
Theorem 3.4) and shows therefore that the number of modes is identifiable from
data, but consistency is not uniform, and for fixed $n$ only a lower bound for
$J(P)$ can be given.

The key ingredient for Donoho's results (as well as the result
of \cite{BahSav56}) is the richness of
${\cal P}$. If ${\cal P}$ is suitably constrained by assumptions, the 
functionals can be empirically identified; however it cannot be empirically
identified
whether such assumptions hold.
\end{example}

\begin{example}\label{excausality} \cite{SpGlSc93,RSSW03} deal with the 
possibility to infer the presence 
or absence of causal arrows in a graphical model. 
Chapter 4 of \cite{SpGlSc93} is about ``statistical indistinguishability''. They
define several indistinguishability concepts of different strengths for the 
problem of identifying the causal graph. As classical identifiability, these
concepts regard the model, not making reference to observable data, but 
empirical identifiability is also of key interest regarding the issue how much 
can be inferred about the causal relationship between observed variables in the
presence of an unobserved confounder.   

\cite{RSSW03} show that in several
situations  ``uniformly consistent tests'' do not exist, whereas  
``pointwise consistent tests'' exist; the latter 
however do not allow to distinguish 
presence or absence of causal arrows for any fixed $n$ uniformly over the
possible parameters. ``Consistent tests'' are procedures that can have outcomes
0 (``accept absence of arrow''), 1 (``reject absence or arrow''), or 2
(``inconclusive''). Their Example 2 is very simple and most instructive. Assume
observed binary random variables $X$ and $Y$ and a categorical unobserved 
confounder $Z$. The existence of a causal arrow between $X$ and $Y$ is 
operationalised by a parameter $\theta^*$ encoding the strength of the causal
effect, where $\theta^*=0$ means that there is no causal arrow.
There are eight different possible causal graphs encoding the possible 
conditional independence structures (Fig. 3 in \cite{RSSW03}). The authors 
assume that the distribution is faithful to the graph, meaning that there are 
not more independence relationships in the distribution than encoded in the 
graph. 

Adapting the
terminology of the present work, the problem is to distinguish two classes of 
graphs. The first class $C_1$ consists of those graphs 
that encode marginal independence between $X$ and $Y$,
which implies $\theta^*=0$ under faithfulness.
The second class $C_2$ 
are the graphs that imply marginal dependence between $X$ and
$Y$, in which case a consistent test should give an inconclusive result, 
because there is a possible graph in which both $X$ and $Y$ are influenced
by $Z$, which causes marginal dependence, despite $X$ and $Y$ being independent
given $Z$, therefore $\theta^*=0$. The results in 
\cite{RSSW03} imply that within the second
class, the existence of a causal arrow between $X$ and $Y$ is not empirically identifiable, as only dependence or independence between $X$ and $Y$ can be observed. There is a pointwise consistent test (i.e., a consistent estimator of the indicator variable, therefore empirical identifiability) that can tell 
apart the first and the second class. The two classes  are not empirically gap 
empirically distinguishable, but this is not very surprising as $\theta$ can 
be arbitrarily close to 0 if a causal arrow exists. What is more remarkable is 
that the proof of their Theorem 1
(which states that no uniformly consistent test exists) implies that $C_1$ 
cannot even be empirically gap distinguished from $C_2\cap \{\theta^*=\theta^*_0\}$, where
$\theta^*_0\neq 0$ is a fixed parameter value. This is 
because a graph in the second class 
that has causal arrows between each pair of $X$, $Y$, and $Z$ encodes a model
that allows for dependence between $Z$ and $X$, 
and also between $Z$ and $Y$ in such a way that $X$ and $Y$ can 
be arbitrarily close to marginally independent 
despite the existing causal effect $\theta^*_0$.
\end{example} 

\begin{example}\label{emissing} 
Regarding models for missing values, a key distinction is 
between MAR (missing at random) and MNAR (missing not at random) mechanisms.
MAR holds if the distribution of the missingness indicator only depends on the 
complete data (including missing values) through the non-missing observations.
This is a very convenient assumption for dealing with missing values, because
it means that the non-missing data provide enough information to allow for 
unbiased inference. As acknowledged by the missing values literature, 
it is doubtful that this assumption is realistic, though,
and it is also doubtful whether the data allow to check this assumption, as
key information for this is hidden in the missing values. 

\cite{MBSK08} made this concern precise by showing that for every MAR model
there exists an MNAR model that reproduces the same observed likelihood function,
meaning in particular that the densities of the observed data and therefore
the probabilities for every observable set are equal between these models.
This translates into empirical indistinguishability (not even potential 
distinguishability is possible) 
of an indicator of MAR vs. MNAR, once more 
assuming a sufficiently rich class of models. The authors state that MAR and
MNAR may be distinguishable under certain parametric assumptions that restrict
the flexibility of the MNAR models to emulate the likelihoods for certain
MAR models; but then it is not empirically distinguishable 
whether such a model holds or not. 
\end{example}
\begin{example}\label{eordinal} $p$-dimensional ordinal data is often modelled as generated by
discretisation of latent Euclidean 
continuous variables. There is much work about 
dimension reduction, but assume here that there are $p$ continuous 
latent variables, every
one of which corresponds to an observed ordinal variable, i.e., if $Y_i$ is the
$i$th latent variable, and the observed ordinal variable $X_i$ takes the
ordered 
categories $j=1,\ldots,k_i$ with probabilities $\pi_1,\ldots,\pi_{k_i}$ with $\pi_0=0$, $X_i=j$ if $Y_i$ is between the $\sum_{l=0}^{j-1}\pi_l$- and 
$\sum_{l=0}^{j}\pi_l$-quantiles of ${\cal L}(Y_i)$. 
  
The most popular approach is to assume the latent variables as multivariate 
Gaussian, see
\cite{Muthen84}. A mis-specification of the distribution of the latent 
distribution can have consequences in practice, see \cite{FolGro20}, who 
discuss tests of this assumption. In fact, all such tests only test the 
dependence structure of $(Y_1,\ldots,Y_p)$ rather than the shape
of the marginal distributions of ${\cal L}(Y_i)$, which makes sense as the 
assignment of categories of $X_i$ according to quantiles of ${\cal L}(Y_i)$ 
obviously works for any continuous ${\cal L}(Y_i)$. Following
Sklar's famous theorem (\cite{Sklar59}), the joint distribution of $(Y_1,\ldots,Y_p)$ has a cumulative distribution function (CDF)
$H$ that can be written as $H(y_1,\ldots,y_p)=C(F_1(y_1),\ldots,F_p(y_p))$ where
$F_1,\ldots,F_p$ are the marginal CDFs and $C$ is a copula. 
This means that
any dependence structure observed in ordinal data is compatible with any choice
of the marginals.
Indeed, for the latent 
variable model for such data, Proposition 2 of \cite{AlmMou14} implies that 
$F_1,\ldots,F_p$ are not empirically identifiable, and that 
any two vectors of marginal CDFs are not even empirically 
distinguishable (be it 
potentially). The authors formulate this as classical identifiability
statement, but involving the $(Y_1,\ldots,Y_p)$ in the model, even though 
not observable, means that the model is identified but empirically, i.e., 
not from what is observable.
\end{example}

\section{Distinguishing independence and dependence}\label{sindep}
The problem of identifying constant correlation between Gaussian observations is
an instance of the more general problem to detect dependence between the 
observations in a sample, particularly if they are meant to be analysed by
methods that assume independent data. Here the focus is on i.i.d. data. There
will not be sophisticated results in this Section, the focus is on general ideas.

Existing tests such as the runs test (\cite{WalWol40})  
require additional information about the kind of dependence to be detected. 
Most of them test for dependence governed by the observation order, which is
sensible if it can be suspected that closeness in observation order can give
rise to dependence. This is often the case if observations are a time series, 
but also other meaningful orderings are conceivable, and 
also dependence governed by 
``closeness''  on external variables such as spatial location. 
Alternatively, there may be a known grouping of
observations and possible within-group dependence, as modelled for example
by random effects models. 
 
In practice, the observation order is not always meaningful, an originally 
existing meaningful observation order may be unavailable to the data analyst, 
or dependence 
structures can be suspected that cannot be detected by examining relations
between observations that are in some sense ``close''. Constant correlation 
between any two Gaussian observations as in model M1
is one example of such a structure.

The question of interest here is whether independence and dependence can be
distinguished in case that the observation order carries no relevant 
information, and
neither is there secondary information from additionally observed variables.
This amounts to observing the 
empirical distribution of the data only. 

For real-valued data $X_1,\ldots,X_n$, $\hat F_n$ denotes the empirical 
distribution function. Assume that only $\hat F_n$ is observed. The concept of
empirical identifiability can be applied to binary ``parameters'', and 
particularly to a parameter that indicates, within an underlying model, whether
$X_1,\ldots,X_n$ are i.i.d. or not. Empirical identifiability involves 
asymptotics. For $n\to\infty$ here it needs to be
assumed that a new sequence $X_1,\ldots,X_n$ is generated for each $n$, because
observing a sequence $\hat F_1,\ldots,\hat F_n,\hat F_{n+1},\ldots$ based on 
the same sequence of observations will re-introduce observation order 
information.   

At first sight the task of identifying dependence from $\hat F_n$ may seem 
hopeless, given that $\hat F_n$ is perfectly compatible with i.i.d. data 
generated from distributions with a ``close'' true CDF or even $\hat F_n$
itself. As in other identifiability problems, information can come from
restrictive assumptions. 
\begin{example}\label{edepmixture}
Consider model M2 and assume for the number of groups $m\ge 2$, but only $\hat F_n$ is observed, meaning that it cannot be observed to which group $i$ an observation $X_{ij}$ belongs. Independence amounts to 
$\tau_1^2=0$, and the interest here is in $g(\theta)=1(\tau_1^2=0)$ denoting the indicator function for $\{\tau_1^2=0\}$. In case that $\tau_1^2=0$, the underlying distribution of $X_{ij}$ is i.i.d. Gaussian. In case that $\tau_1^2>0$, the underlying distribution of $X_{ij}$ partitions the data into different Gaussians for different $i$. Assuming that $\frac{n_i}{n}\to \pi_i>0$, $\hat F_n$ will for large enough $n$ look like a Gaussian mixture, which can be told apart from a single Gaussian. 
The order of a Gaussian mixture can be consistently estimated (\cite{JaMaPr01}), and therefore $1(\tau_1^2=0)$, which is equal to the indicator of a single mixture component, is empirically identifiable. Note that the cited result is for i.i.d. data from a Gaussian mixture, which in model M2 would require the group memberships to be modelled i.i.d. multinomial$(1,\pi_1,\ldots,\pi_m)$, and then conditioning on the unknown values of $Z_1,\ldots, Z_m$.

The example illustrates that certain empirical distributions can indeed indicate dependence, if corresponding distributional shapes (here a Gaussian mixture with $m\ge 2$ components) are assumed as impossible under independence but can occur under dependence.
\end{example}   
Even if general marginal distributions are allowed, there are specific dependence structures that can be identified from the empirical distribution alone. In order to simplify matters, from now on consider binary data $X_1,\ldots,X_n$, for which observing the empirical distribution is equivalent to observing the number of ones or $\bar X_n=1-\hat F_n(0)$, and all marginal distributions are Bernoulli$(p)$. Call the i.i.d. model M3. A problem of interest is whether there are models for dependence for which all marginal 
distributions are identical that can be told apart based on $\hat F_n$ 
from M3. 
\begin{example}
Here is an example for a dependence structure that can be distinguished from $\hat F_n$. This relies on $\bar X_n$ concentrating on a pre-speficied value (here $\frac{1}{2}$) with larger probability than under M3 even if $p$ equals this value. Consider a model M4 where ${\cal L}(X_1)=$Bernoulli$(p)$. With $q=p1\left(p\le \frac{1}{2}\right)+(1-p)1\left(p> \frac{1}{2}\right)$, 
$r=2\frac{q}{q+1},$ let $Y$ be an unobserved Bernoulli$(r)$-random variable. 
If $Y=0$, $X_2,\ldots,X_n,\ldots$ i.i.d. Bernoulli$(p/2)$.
If $Y=1$, $X_3,X_5, X_7,\ldots$ i.i.d. Bernoulli$(1/2)$, and for all even 
$n:\ X_n=1-X_{n-1}$, so that $\bar X_n=\frac{1}{2}$. $r$ is chosen
so that all marginal distributions are Bernoulli$(p)$. With
$A_n=\{\bar X_n=\frac{1}{2}\}$, $P(A_n)\to 0$ under M3, whereas for even $n$,
$P(A_n)=r$ under M4. Therefore $A_n$ gap distinguishes M3 from M4.
\end{example}
Such examples rely on the definition of a subset $A_n$ of possible
values of $\hat F_n$, which under the dependence model for large enough $n$
has a probability either higher than the maximum over $p$ under M3, or smaller than the corresponding minimum. For checking dependence in practice, $A_n$ needs to be specified in advance, meaning that the user needs to know a priori which values of $\hat F_n$ can be suspected to indicate dependence even given the possibility under M3 that $p=\bar X_n$. Such information is rarely available in practice. 

Summarising, dependence can be diagnosed from the data only if 
\begin{itemize}
\item it is governed by the known order of observations or known external variables, 
\item or it favours (or avoids) specific events regarding the observed empirical distribution compared to the independence model of interest, 
\end{itemize}
both in ways that the user has to specify in advance. It can be suspected that many existing dependence structures are not of this kind, meaning that only very limited aspects of independence between observations, regarding a sequence of such observations, can be checked. It is therefore very important to think through all background information about the data generating process to become aware of further potential issues with independence. 

\section{Cluster membership in $k$-means clustering}\label{skmeans}
$k$-means clustering is probably the most popular cluster analysis method 
(\cite{Jain10}). It can be connected to a ``fixed classification model'' 
(\cite{Bock96}):
Let $\X_1,\ldots,\X_n,\ \X_i\in\R^p,\ i=1,\ldots,n,$ be 
independently distributed with 
\begin{equation}\label{efpnormal}
{\cal L}(\X_i)={\cal N}_p\left(\bmu_{\gamma_i},\sigma^2\I_p\right),\ \gamma_i\in\{1,\ldots,k\},\ k>1,\ \sigma^2\ge 0. 
\end{equation}
This model can be interpreted as generating $k$ different Gaussian distributed clusters characterised by cluster means $\bmu_1,\ldots,\bmu_k\in\R^p$, all with the same spherical covariance matrix, and $\gamma_i$ indicates the true cluster membership of $\X_i$. The $\gamma_i$ take discrete values, and their number converges to $\infty$ with $n$, so these are nonstandard parameters, but in many applications they are of practical interest. 

The maximum likelihood (ML)-estimator for $\theta=(\bmu_1,\ldots,\bmu_k,\gamma_i,\ldots,\gamma_n)$ in this model is given by $k$-means clustering (the ML-estimator for $\sigma^2$ is easily derived, but this is not relevant here) of data $\tilde{\X}_n=(\X_1,\ldots,\X_n)$:
\begin{eqnarray*}
 T_n(\tilde{\X}_n)&=&(\m_{1n},\ldots,\m_{kn},g_{in},\ldots,g_{nn})\\
&=&
\argmin_{\m_1,\ldots,\m_k,g_1,\ldots,g_k} W(\m_{1n},\ldots,\m_{kn},g_{in},\ldots,g_{nn}),\\
 W(\m_{1n},\ldots,\m_{kn},g_{in},\ldots,g_{nn})&=&  
\sum_{i=1}^n \|\X_{in}-\m_{g_{in}}\|^2
\end{eqnarray*}
with ties in the $\argmin$ broken in an arbitrary way.
For given $\m_{1},\ldots,\m_{k}$, the $g_1,\ldots,g_{n}$ minimising $W$ 
are given by 
\[
g_i=\argmin_{j\in\{1,\ldots,k\}}\|\X_i-\m_j\|,\ i=1,\ldots,n,
\] 
and with these 
write $W(\m_{1n},\ldots,\m_{kn})=W(\m_{1n},\ldots,\m_{kn},g_{in},\ldots,g_{nn})$.
For given $g_1,\ldots,g_{n}$, the cluster-wise mean vectors minimise $W$. As there are only finitely many values of $g_1,\ldots,g_{n}$, the 
ML-estimator does always exist, if not necessarily uniquely. 
Two issues with identifiability here are (i) that the numbering of the clusters is arbitrary and (ii) that $\bmu_q=\bmu_r$ for $q\neq r$ means that it is not possible to distinguish between $\gamma_i=q$ and $\gamma_i=r$ for $i\in\N$. Therefore assume that the $\bmu_1,\ldots,\bmu_k$ are pairwise different and lexicographically ordered (i.e., with obvious notation, $\mu_{11}\le\ldots\le\mu_{k1}$ with ties broken by the second variable, or, if there's still equality, by the third and so on, same for $\m_{1n},\ldots,\m_{kn}$). This makes the model identifiable according to the traditional definition.  

However, due to the nonstandard nature of the model parameters, the ML-estimator is known to be inconsistent, even if only the estimation of the $k$ mean vectors alone is of interest (\cite{Bryant91}). 

The cluster membership parameters are another example for parameters that are identifiable according to the classical definition (because $\gamma_i$ uniquely defines the distribution of $\X_i$), but cannot be empirically identified.
\begin{theorem}\label{tgaussfpident} The parameters $\gamma_i,\ i\in\N$ in the fixed classification model defined in (\ref{efpnormal}) are not empirically identifiable. 
\end{theorem}

It may be suspected that this is a consequence of the fact that for $i=1,\ldots,n$, only $\X_i$ holds information about the parameter $\gamma_i$, and the number of these parameters goes to $\infty$ with $n\to\infty$. But this is not quite true. More observations add information about the clusters that can in turn be used to classify individual observations better. The problem here is rather the Gaussian distribution assumption, which implies that the marginal density of $\X_i$ is everywhere nonzero, so that the single observation made of $\X_i$ is not enough to determine with probability 1 to what cluster the observation belongs (the setup by which the classical identifiability could be used to estimate this parameter would be to have a potentially infinite amount of replicates of $\X_i$), even if there is an infinite amount of information about the clusters. In fact, there is a different model setup in which the $\gamma_i$ are empirically identifiable, which requires that, where densities exist, 
the marginal density $f_{\theta^*,n}(\X_i=\x)$ is zero wherever $f_{\theta,n}(\X_i=\x)>0$, where $\theta^*$ equals $\theta$ in all components except $\gamma_i$.

Defining 
\[
W(P)=(\bmu_1^*,\ldots,\bmu_k^*)=\argmin_{(\m_1,\ldots,\m_k)\in (\R^p)^k} \int \min_{\m\in\{\m_1,\ldots,\m_k\}}\|\x-\m\|^2 dP(\x),
\] 
\cite{Pollard81} showed that for a distribution $P$ satisfying 
\begin{equation}\label{epollardass}
E_P\|\X\|^2 < \infty,\ W(P) \mbox{ is unique up to the numbering of the means,}
\end{equation}
the $k$-means estimator $\left(T_n^m\right)_{n\in\N}$, where $T_n^m(\tilde{\X}_n)=(\m_{1n},\ldots,\m_{kn})$, is strongly consistent for $W(P)$. Assume further that 
\begin{equation}\label{econtass}
\forall j\neq l\in\{1,\ldots,k\}:\ P\{\|\X-\bmu^*_j\|^2=\|\X-\bmu^*_l\|^2\}=0.
\end{equation}
For ${\cal L}(\X)=P,\ j=1,\ldots,k$, define 
\[
A_j=\left\{\X:\ j=\argmin_l\|\X-\bmu_l^*\|^2\right\},\ 
P_j={\cal L}(\X\mid\X\in A_j),\ \pi_j=P(A_j).
\]
$P_j$ is $P$ constrained to the set $A_j$ of points that are closest 
to the mean 
$\bmu_j$ ($A_1,\ldots,A_k$ form a Voronoi tesselation of $\R^p$), and
\begin{equation}\label{epmix}
P=\sum_{j=1}^k\pi_jP_j
\end{equation}
 (every distribution can be written as a mixture in 
this form; as a side remark, $P$ might be a Gaussian mixture, but in this case 
the $P_j$ are not its Gaussian components). Mixture models of this form 
can be derived from a model for outcomes $(G, \X)$ with $G\in\{1,\ldots,k\}$
distributed according to a categorical distribution with probabilities
$(\pi_1,\ldots,\pi_k)$ and ${\cal L}(\X\mid G=j)=P_j$. Then ${\cal L}(\X)=P$ 
(\cite{McLPee00}). For an i.i.d. sequence $\Y=((G_1,\X_1),(G_2,\X_2),\ldots)$
let ${\cal L}(\Y)=\tilde P,\ \G=(G_1,G_2,\ldots)$.

Now consider ${\cal L}(\tilde{\X}_n)=P^*$ so that $\X_1,\ldots,\X_n$  are
independently distributed with 
\begin{equation}\label{efpgeneral}
{\cal L}(\X_i)=P_{\gamma_i}\ \gamma_i\in\{1,\ldots,k\},\ k>1,\ i=1,\ldots,n. 
\end{equation}
This defines a fixed classification model associated to the mixture $P$.  
Let $Q_P$ be an infinite i.i.d. product of categorical distributions on 
$\{1,\ldots,k\}$ with probabilities $(\pi_1,\ldots,\pi_k)$. Assume for given 
$P$ that $\bm{\gamma}=(\gamma_1,\gamma_2,\ldots)$ fulfill 
\begin{equation}\label{econslabels}
\tilde P\left\{\lim_{n\to\infty}T_n^m(\tilde{\X}_n)=W(P)\mid \G=\bm{\gamma}\right\}=1.
\end{equation}
Because of the strong consistency of $T_n^m$, (\ref{econslabels}) holds with
probability 1 under $Q_P$, but note that under (\ref{efpgeneral}), 
$\bm{\gamma}$ is a fixed 
parameter and not a random variable, and the fact that (\ref{econslabels}) 
holds for $Q_P$-almost all $\bm{\gamma}$ just means that (\ref{econslabels}) is not 
more restrictive than assuming a mixture with fixed proportions, although it 
will not allow for fully general $\bm{\gamma}$.
\begin{theorem}\label{tidentpollard}  Assuming  
(\ref{epollardass}), (\ref{econtass}), and (\ref{econslabels}),
the parameters $\gamma_i,\ i\in \N,$ in the fixed classification model defined 
by (\ref{efpgeneral}) are empirically identifiable. 
\end{theorem}
Already from \cite{Pollard81} it is clear that $k$-means
does not actually estimate
the centres of the spherical Gaussians in 
(\ref{efpnormal}), but rather the Voronoi tesselation resulting from $P$, and
the resulting clusters are not necessarily spherical. Added here is the
observation that one can define meaningful cluster indicators in this setup,
and that these can be consistently estimated, even though there is one such
indicator for every observation. This is not possible in 
(\ref{efpnormal}).
Furthermore 
(\ref{epmix}) interprets $P$ as a mixture,
and thus shows that there is a mixture
that $k$-means estimates consistently.

The reader may wonder about empirical distinguishability of the parameters $\gamma_i,\ i\in\N$, i.e., about whether the given values $j_1$ and $j_2$ of $\gamma_i$ for given $i$ could be distinguished. In the situation of Theorem \ref{tidentpollard}, this follows from Lemma \ref{lindist} (projections of discrete parameters are by definition uniformly continuous). The situation of Theorem \ref{tgaussfpident} is less obvious. It is however clear that the data contain some information about $\gamma_i$ through $\X_i$. If it were possible to empirically identify $\bmu_1,\ldots,\bmu_k$, $A=\{\|\X_i-\bmu_{j_1}\|>\|\X_i-\bmu_{j_2}\|\}$ could distinguish $j_1$ and $j_2$. A conjecture is that finding consistent estimators for $\bmu_1,\ldots,\bmu_k$ requires additional conditions on the sequence  $(\gamma_i)_{i\in\N}$ that allow to use a consistent estimator from an i.i.d. mixture, see \cite{RedWal84,Bryant91}.

\section{Conclusion} \label{sconc}
There are various potential issues with the identification of parameters, and the four definitions given here (empirical identifiability, empirical distinguishability, empirical gap distinguishability, potential empirical distinguishability) may not cover all of them; even using the definition of distinguishing sets, further definitions are possible, for example empirical set identifiability, but what is already present allows to deal with many examples. 

Apart from the precise definitions, there are also different sources for identifiability and distinguishability problems. In some situations (Examples \ref{emissing}, \ref{eordinal}) the problem is that modelled information is not observed, either because it regards missing values, or latent variables. In some situations (Examples \ref{ebahsav}, \ref{edensity}), the issue is that the class of distributions to consider, even for a single parameter value of interest, is too rich, allowing for so much flexibility that the probability of any observable set cannot be sufficiently constrained. Identifiability and empirical distinguishability can be the result of model assumptions constraining this flexibility, see Example \ref{edepmixture}; in Example \ref{exgausscor}, constraining the variance $\sigma^2$ would in turn allow the data to be informative about $\rho$. These model assumptions cannot be justified from the data alone though. In some further situations, the data carry information about the parameter of interest, but this information, or more precisely the growth of the information over $n$, is limited (Lemma \ref{lrandomeff}, Example \ref{eneysco}, Theorem \ref{tgaussfpident}; Example \ref{exgausscor} is also of this kind, using (\ref{erandomeff})). Example \ref{exindist} is constructed so that the distinguishing information may not occur at any finite $n$, and the parameterisation in Example \ref{excausality} is arbitrarily close to a situation of classical non-identifiability, which is only avoided by the faithfulness assumption.  

Some examples such as model M01 and the problem in Section \ref{skmeans} are characterised by not allowing for i.i.d. repetition; the corresponding parameters can be identified if the whole sequence of observations is repeated i.i.d., and the lack of empirical identifiability is due to the assumed impossibility to do this. It may be wondered whether such models have a valid frequentist interpretation, which seems to rely on i.i.d. replicability at least in principle. Frequentism needs to be interpreted in a rather ``idealist'' way to accommodate such situations, appealing to replication of an infinite sequence as a thought experiment, although this issue can arguably be made regarding time series and other models as well; for more on this, see \cite{Hennig20}.   

The most relevant and unsettling implication of the work for practice regards the lack of possibility to check certain model assumptions, particularly independence; flexible enough models allowing for non-identical marginals can be impossible to detect as well, although this is not shown here.

The considerations regarding requirements for detecting dependence in Section \ref{sindep} do not only hold for data for which only the empirical distribution is observed; they hold in the same way also for situations in which the observation order, even if known and potentially meaningful, is not informative for the dependence structure, and no external variables exist either that hold such information. This is probably a very common situation. The only way to justify independence then is knowledge about the subject matter and the data generation. Bayesians may think that the lack of information in the data about parameters such as $\rho$ in M01 could be compensated by a prior distribution, but the lack of empirical identifiability and distinguishability of $\rho$ raises the question where quantitative information should come from to set up the prior. A prior could only be obtained from existing qualitative information about the data generating process, and as there is no information in the data, the prior will determine the impact of $\rho$ without the possibility of being ``corrected'' by the data.

\section*{Appendix: Proofs}
\subsection*{Proof of Lemma \ref{lcondgauss}}
Recall that, in general, if ${\cal L}\begin{pmatrix}\Z_1\\ \Z_2\end{pmatrix}={\cal N}\left(\begin{pmatrix}\bmu_1\\\bmu_2\end{pmatrix},
\begin{bmatrix}\bS_1&\bS_{12}\\ \bS_{12} & \bS_{2}\end{bmatrix}\right)$, then
\begin{equation}\label{econdnorm}
{\cal L}(\Z_1\mid \Z_2)={\cal N}\left(\bmu_1+\bS_{12}\bS_{2}^{-1}(\Z_2-\bmu_2),\bS_1-\bS_{12}\bS_{2}^{-1}\bS_{12}\right).
\end{equation}
Now consider 
\[
{\cal L}\begin{pmatrix} X_1\\ \vdots \\ X_n\\\bar{X}_n\end{pmatrix}=
{\cal N}_{n+1} \left(\begin{pmatrix}\mu\\ \vdots\\ \mu\end{pmatrix},
\begin{bmatrix}
\sigma^2 & \rho\sigma^2 & \dots & \rho\sigma^2 & \tilde{\sigma}^2\\
    \rho\sigma^2 & \sigma^2 & \dots & \rho\sigma^2 & \tilde{\sigma}^2\\
    \vdots & \vdots & \ddots &\vdots & \vdots\\
    \rho\sigma^2 & \rho\sigma^2 & \dots & \sigma^2 & \tilde{\sigma}^2\\
    \tilde{\sigma}^2 & \tilde{\sigma}^2 & \dots
    & \tilde{\sigma}^2 & \tilde{\sigma}^2
    \end{bmatrix}\right).
\]
Let $\tilde{\sigma}^2=\frac{\sigma^2(1+(n-1)\rho)}{n}>0$.
Using (\ref{econdnorm}), ${\cal L}(\X_n\mid\bar{X}_n)={\cal N}_n(\tilde{\bmu},\tilde{\bS})$, where 
\begin{eqnarray*}
\tilde{\bmu} &=& \bmu+\begin{pmatrix} \tilde{\sigma}^2\\ \vdots\\ \\ \tilde{\sigma}^2\end{pmatrix}\frac{1}{\tilde{\sigma}^2}(\bar{X}_n-\mu)=\begin{pmatrix} \bar{X}_n\\ \vdots\\ \bar{X}_n\end{pmatrix},\\
\tilde{\bS} &=& \bS-\begin{pmatrix} \tilde{\sigma}^2\\ \vdots\\ \\ \tilde{\sigma}^2\end{pmatrix}\frac{1}{\tilde{\sigma}^2}(\tilde{\sigma}^2,\ldots,\tilde{\sigma}^2)\\
&=& \bS-\begin{bmatrix} \tilde{\sigma}^2 & \dots & \tilde{\sigma}^2\\
\vdots & \ddots & \vdots\\ \tilde{\sigma}^2 & \dots & \tilde{\sigma}^2
\end{bmatrix}= \begin{bmatrix}\left(1-\frac{1}{n}\right)\sigma^2_* & -\frac{1}{n}\sigma^2_* & \dots &  -\frac{1}{n}\sigma^2_*\\
-\frac{1}{n}\sigma^2_* & \left(1-\frac{1}{n}\right)\sigma^2_* & \dots & -\frac{1}{n}\sigma^2_*\\
\vdots & \vdots & \ddots & \vdots\\
-\frac{1}{n}\sigma^2_* & -\frac{1}{n}\sigma^2_* & \dots & \left(1-\frac{1}{n}\right)\sigma^2_* \end{bmatrix}
\end{eqnarray*}
with $\sigma^2_*=(1-\rho)\sigma^2$. Applying (\ref{econdnorm}) once more, ${\cal L}(\Y_n\mid\bar{Y}_n)={\cal N}_n(\bmu_*,\bS_*)$, where
\begin{eqnarray*}
\bmu_* &=& \bmu+\begin{pmatrix} \frac{\sigma_*^2}{n}\\ \vdots\\ \\ \frac{\sigma_*^2}{n}\end{pmatrix}\frac{n}{\sigma_*^2}(\bar{Y}_n-\mu)=\begin{pmatrix} \bar{Y}_n\\ \vdots\\ \bar{Y}_n\end{pmatrix},\\
\bS_* &=& \sigma_*^2{\bf I}_n-\begin{pmatrix} \frac{\sigma_*^2}{n}\\ \vdots\\ \\ \frac{\sigma_*^2}{n}\end{pmatrix}\frac{n}{\sigma_*^2}\left(\frac{\sigma_*^2}{n},\ldots,\frac{\sigma_*^2}{n}\right)\\
&=& \begin{bmatrix}\left(1-\frac{1}{n}\right)\sigma_*^2 & -\frac{1}{n}\sigma_*^2 & \dots &  -\frac{1}{n}\sigma_*^2\\
-\frac{1}{n}\sigma_*^2 & \left(1-\frac{1}{n}\right)\sigma_*^2 & \dots & -\frac{1}{n}\sigma_*^2\\
\vdots & \vdots & \ddots & \vdots\\
-\frac{1}{n}\sigma_*^2 & -\frac{1}{n}\sigma_*^2 & \dots & \left(1-\frac{1}{n}\right)\sigma_*^2 \end{bmatrix}.\\
\mbox{Therefore, }&&{\cal L}(\X_n\mid\bar{X}_n)={\cal L}(\Y_n\mid\bar{Y}_n),
\end{eqnarray*} 
and these do not depend on $\mu$.
\subsection*{Proof of Theorem \ref{tnotidrho}}
Let $(R_n)_{n\in\N}$ with $R_n:\ {\cal X}^n\mapsto [0,1]$ 
be a consistent estimator of $\rho$ ($\rho=1$ and $R_n=1$ are assumed 
possible, although the following contradiction to consistency relies on 
$\rho<1$; $R_n$ could be allowed to take negative values even though 
$\rho\ge 0$ must be assumed because otherwise the matrices $\bS_j$ below would 
become negative definite for large enough $n$).

For $j=1,2,$ consider ${\cal L}(\X_n^{(j)})= {\cal N}_n(\bm{\mu,\Sigma}_j)$, 
where
$\bS_1$ defined by $0\le\rho=\rho_1<1,\ \sigma^2=\sigma_1^2>0$ and $\bS_2$ defined by $\rho_1<\rho=\rho_2<1,\ \sigma^2=\sigma_2^2=\frac{1-\rho_1}{1-\rho_2}\sigma_1^2.$ With ${\cal L}(\Y_n)={\cal N}_n(\bmu,(1-\rho_1)\sigma_1^2\I_n)$, Lemma \ref{lcondgauss} yields
\begin{equation}\label{ethreeequal} 
{\cal L}(\X_n^{(1)}\mid\bar{X}_n^{(1)})={\cal L}(\Y_n\mid\bar{Y}_n)={\cal L}(\X_n^{(2)}\mid\bar{X}_n^{(2)}).
\end{equation}
Let $\epsilon<(\rho_2-\rho_1)/2$, and $A=[\rho_1-\epsilon,\rho_1+\epsilon]$. Consistency implies
\begin{equation}\label{econsistency}
\lim_{n\to\infty} P^{\X^{(1)}_n}\{R_n\in A\}=1,\mbox{ but }
\lim_{n\to\infty}P^{\X^{(2)}_n}\{R_n\in A\}=0.
\end{equation}
This implies that there must be a set $B\subseteq \R$ so that 
\[ 
\lim_{n\to\infty} P^{\bar{X}^{(1)}_n}(B)=1,\ \forall x\in B:\ \lim_{n\to\infty} P^{\X^{(1)}_n}\{R_n\in A\mid\bar{X}^{(1)}_n=x\}=1.
\]
The conditional distribution ${\cal L}(\X_n^{(1)}\mid\bar{X}_n^{(1)}=1)$ is determined up to a null set of values of $\bar{X}_n^{(1)}$ according to the disintegration theorem. This means that the required $B$ exists 
for any valid choice of ${\cal L}(\X_n^{(1)}\mid\bar{X}_n^{(1)}=x)$ for all $x$. The same argument holds later when conditioning on sets with probability zero.

By (\ref{ethreeequal}) and (\ref{econsistency}):
\[ 
\forall x\in B:\ 
\lim_{n\to\infty} P^{\X^{(2)}_n}\{R_n\in A\mid\bar{X}^{(2)}_n=x\}=1 \Rightarrow
\lim_{n\to\infty} P^{\bar{X}^{(2)}_n}(B)=0.
\]
Now 
\begin{eqnarray*}
{\cal L}(\bar{X}_n^{(1)})&=&{\cal N}\left(\mu,\frac{(1-\rho_1)\sigma^2_1}{n}+\rho_1\sigma^2_1\right)\to_{n\to\infty} {\cal N}\left(\mu,\rho_1\sigma^2_1\right),\\
{\cal L}(\bar{X}_n^{(2)})&=&{\cal N}\left(\mu,\frac{(1-\rho_2)\sigma^2_2}{n}+\rho_2\sigma_2^2\right) \to_{n\to\infty} {\cal N}\left(\mu,\frac{\rho_2(1-\rho_1)}{1-\rho_2}\sigma^2_1\right),
\end{eqnarray*}
both with non-vanishing variance, and therefore $B$ with 
\[
\lim_{n\to\infty} P^{\bar{X}^{(1)}_n}(B)=1,\ \lim_{n\to\infty} P^{\bar{X}^{(2)}_n}(B)=0
\]
 cannot exist, with contradiction, so $R_n$ cannot be consistent.
\subsection*{Proof of Theorem \ref{tmunonid}}
Let $(T_n)_{n\in\N}$ be a consistent 
estimator of 
$\mu$. 
For $j=1,2,$ consider ${\cal L}(\X_n^{(j)})= {\cal N}_n(\bmu_j,\bS)$.
For $\mu_1\neq\mu_2$, let $\epsilon<(\mu_2-\mu_1)/2$, and $A=[\mu_1-\epsilon,\mu_1+\epsilon]$. Consistency implies
\begin{equation}\label{econsistencymu}
\lim_{n\to\infty} P^{\X^{(1)}_n}\{T_n\in A\}=1,\mbox{ but }
\lim_{n\to\infty}P^{\X^{(2)}_n}\{T_n\in A\}=0.
\end{equation}
Now consider the equivalent model formulation (\ref{erandomeff}), for $j=1, 2,$
denoting the $Z$ as $Z^{(j)}$.
There must be a set $B\subseteq \R$ so that 
\[ 
\lim_{n\to\infty} P^{Z^{(1)}}(B)=1,\ \forall z\in B:\ \lim_{n\to\infty} P^{\X^{(1)}_n}\{T_n\in A\mid Z^{(1)}=z\}=1.
\]
By (\ref{erandomeff}), with $\z=(z,\ldots,z)'$ and 
$\bmu_1=(\mu_1,\ldots,\mu_1)'$,
\[
{\cal L}(\X_n^{(1)}\mid Z^{(1)}=z)={\cal N}(\bmu_1+\z,\bS)=
{\cal L}(\X_n^{(2)}\mid Z^{(2)}=z+\mu_1-\mu_2).
\]
Therefore, from (\ref{econsistencymu})
\begin{eqnarray*}
\forall z\in B:\ 
\lim_{n\to\infty} P^{\X^{(2)}_n}\{T_n\in A\mid Z^{(2)}=z+\mu_1-\mu_2\}&=&1 \Rightarrow \\ 
\lim_{n\to\infty} P^{\X^{(2)}_n}\{Z^{(2)}-\mu_1+\mu_2\in B\}&=& 0.
\end{eqnarray*}
But this is not possible, because ${\cal L}(Z^{(2)})$ is just a mean 
shifted version
of ${\cal L}(Z^{(1)})={\cal N}(0,\tau_1^2)$ regardless of $n$, 
and will therefore assign nonzero
probability to everything that has nonzero
probability under ${\cal L}(Z^{(1)})$.
\subsection{Proof of Lemma \ref{lrandomeff}}
Regarding the empirical non-identifiability of 
$g_1(\theta)=\mu$, the proof is basically the same as 
the proof of Theorem \ref{tmunonid}, with the following changes. 
For $k=1,2,$ let ${\cal L}(\X_n^{(k)})$ be defined as in model M2
with $\mu=\mu_k$.
Denote the corresponding $\Z=(Z_1,\ldots,Z_m)$ as $\Z^{(k)}$. For a consistent
$(T_n)_{n\in\N}$,
there must be a set $B\subseteq \R^m$ so that 
\[ 
\lim_{n\to\infty} P^{\Z^{(1)}}(B)=1,\ \forall \z\in B:\ \lim_{n\to\infty} P^{\X^{(1)}_n}\{T_n\in A\mid\Z^{(1)}=\z\}=1,
\]
with $A$ as in the proof of Theorem \ref{tmunonid}.
With $\z=(z_1,\ldots,z_m)'$ and $\bmu_k=(\mu_k,\ldots,\mu_k)',\ k=1,2$,
obtain for $i=1,\ldots,m,\ j=1,\ldots,n_i$:
\[
{\cal L}(X_{ij}^{(1)}\mid\Z^{(1)}=\z)={\cal N}(\mu_1+z_j,\tau_2^2)=
{\cal L}(X_{ij}^{(2)}\mid\Z^{(2)}=\z+\bmu_1-\bmu_2).
\]
The same argument as in the proof of Theorem \ref{tmunonid} 
then establishes that $(T_n)_{n\in\N}$ cannot be consistent.

Obviously no consistent estimator for $g_2(\theta)=\tau_1$ 
can exist because only $m<\infty$ instances of $Z_i$ are generated. 

Regarding $g_3(\theta)=\tau_2^2$, for $n\to\infty$ there must be 
$i_0\in\{1,\ldots,m\}$ with $n_{i_0}\to\infty$. 
For $j=1,\ldots,n_{i_0}$,
${\cal L}(\X_{i_0j}\mid Z_{i_0}=z)={\cal N}(\mu+z,\tau_2^2)$ i.i.d., and a standard
consistent estimator for the Gaussian variance on these observations will  
estimate $\tau_2^2$ consistently.
\subsection*{Proof of Lemma \ref{lindist}}
Let $\epsilon=d_\Lambda(\lambda_1,\lambda_2)/3$.
Uniform continuity of $g$ implies that there is $\delta>0$ so that for $i=1, 2$:
$$
\theta\in B_\delta(\theta_i) \Rightarrow g(\theta)\in B_\epsilon(g(\theta_i)).
$$
Consistency of $T_n$ means that for large enough $n$ and any $\theta_i$ with
$g(\theta_i)=\lambda_i$,
$P_{\theta_i}\left(T_n\in B_\delta(\theta_i)\right)$ is arbitrarily large, 
as is, because of uniform consistency,
$P_{\theta_i}\left(g(T_n)\in B_\epsilon(\lambda_i)\right)$. 
As this also implies that $P_{\theta_1}\left(g(T_n)\in B_\epsilon(\lambda_2)\right)$ is arbitrarily low, $A=\{g(T_n)\in B_\epsilon(\lambda_2)\}$ can distinguish 
$\lambda_1$ and $\lambda_2$.
\subsection*{Proof of Theorem \ref{tgaurhoindist}}
First consider $\rho_1\neq \rho_2$.
 Denote ${\cal L}(X_1,\ldots,X_n)$ in model M01 by 
$P_{n;\mu,\sigma^2,\rho}$. Consider parameters $\mu_1,\ \sigma_1^2$, and any set $A$ so that for given $n:\ P_{n;\mu_1,\sigma_1^2,\rho_1}(A)< P_{n;\mu_1,\sigma_1^2,\rho_2}(A)$ and $P_{n;\mu_1,\sigma_1^2,\rho_2}(A)>\alpha$. 
Show that there exist $\mu_2,\sigma_2^2$ so that
$P_{n;\mu_2,\sigma_1^2,\rho_1}(A)\ge P_{n;\mu_2,\sigma_2^2,\rho_2}(A)$, which means that $A$ cannot distinguish $\rho_1$ and $\rho_2$, and 
it cannot even potentially distinguish $\rho_1$ from $\rho_2$. 

Choose $\sigma_2^2=\frac{1-\rho_1}{1-\rho_2}\sigma_1^2$ so that, according to Lemma \ref{lcondgauss}, for all $\mu:\ P_{n;\mu,\sigma_1^2,\rho_1}(\bullet\mid\bar{X}_n)= P_{n;\mu,\sigma_2^2,\rho_2}(\bullet\mid\bar{X}_n)$. Then,
for any $\mu$: 
\begin{eqnarray*}
P_{n;\mu,\sigma_1^2,\rho_1}(A)&=&\int P_{n;\mu,\sigma_1^2,\rho_1}(A\mid\bar{X}_n=x) dP_{n;\mu,\sigma_1^2,\rho_1}^{\bar{X}_n}(x)\\
&=& 
\int P_{n;\mu,\sigma_2^2,\rho_2}(A\mid\bar{X}_n=x) dP_{n;\mu,\sigma_1^2,\rho_1}^{\bar{X}_n}(x),\\
P_{n;\mu,\sigma_2^2,\rho_2}(A)&=&\int P_{n;\mu,\sigma_2^2,\rho_2}(A\mid\bar{X}_n=x) dP_{n;\mu,\sigma_2^2,\rho_2}^{\bar{X}_n}(x).
\end{eqnarray*}
Recall $P_{n;\mu,\sigma^2,\rho}^{\bar{X}_n}={\cal N}\left(\mu,\frac{(1-\rho)\sigma^2}{n}+\rho\sigma^2\right)$. According to Lemma \ref{lcondgauss}, $P_{n;\mu,\sigma^2,\rho}(\bullet\mid\bar{X}_n=x)$ does not depend on $\mu$. With 
\[
h(x)= P_{n;\mu,\sigma_2^2,\rho_2}(A\mid\bar{X}_n=x),\ \tau_1^2=\frac{(1-\rho_1)\sigma_1^2}{n}+\rho_1\sigma_1^2,\ \tau_2^2=\frac{(1-\rho_2)\sigma_2^2}{n}+\rho_2\sigma_2^2:
\]
\[
P_{n;\mu,\sigma_1^2,\rho_1}(A)=\int h(x)d{\cal N}(\mu,\tau_1^2)(x),\
P_{n;\mu,\sigma_2^2,\rho_2}(A)=\int h(x)d{\cal N}(\mu,\tau_2^2)(x).
\]
The Gaussian density $\varphi_{\mu,\sigma^2}(x)$ with mean $\mu$ and variance $\sigma^2$ is symmetric in $\mu$ and $x$, so that for $x\in\R:$
\begin{eqnarray*}
\int \varphi_{\mu,\sigma^2}(x) d\mu=\int \varphi_{\mu,\sigma^2}(x) dx&=&1,\mbox{ and,}\\
\mbox{for }\tau_1,\ \tau_2:\ 
\int \left(\varphi_{\mu,\tau_1^2}(x)-\varphi_{\mu,\tau_2^2}(x)\right)d\mu &=& 0.
\end{eqnarray*}
Now, 
\[
P_{n;\mu,\sigma_1^2,\rho_1}(A)-P_{n;\mu,\sigma_2^2,\rho_2}(A)=
\int h(x)\left(\varphi_{\mu,\tau_1^2}(x)-\varphi_{\mu,\tau_2^2}(x)\right)dx.
\]
As $0\le h(x)\le 1$, this is integrable over $\mu$, so
\begin{eqnarray}
\int \int h(x)\left(\varphi_{\mu,\tau_1^2}(x)-\varphi_{\mu,\tau_2^2}(x)\right)dxd\mu&=&\nonumber\\
\int h(x)\int\left(\varphi_{\mu,\tau_1^2}(x)-\varphi_{\mu,\tau_2^2}(x)\right)d\mu dx&=&0.\label{eintzero}
\end{eqnarray}
But this means that if $P_{n;\mu_1,\sigma_1^2,\rho_1}(A)< P_{n;\mu_1,\sigma_2^2,\rho_2}(A)$, then there must be a $\mu_2$ so that 
$P_{n;\mu_2,\sigma_1^2,\rho_1}(A)> P_{n;\mu_2,\sigma_2^2,\rho_2}(A)$, and $A$ cannot distinguish $\rho_1$ and $\rho_2$. 
~\\~\\
$\mu_1\neq \mu_2$ can be distinguished by the set 
$A=\{\lvert\bar{X}_n-\mu_1\rvert>\lvert\bar{X}_n-\mu_2\rvert\}$, as  
\[
\forall n,\ \psi=(\sigma^2,\rho):\ P_{n;(\mu_1,\psi)}(A)<\frac{1}{2},\
P_{n;(\mu_2,\psi)}(A)>\frac{1}{2}.
\]
\subsection*{Proof of Theorem \ref{tgaussfpident}}
Consider for $n\to\infty$ a true parameter vector
$$ 
\theta=(\bmu_1,\ldots,\bmu_k,\sigma^2,\gamma_1,\gamma_2,\ldots),\ \sigma^2>0.
$$ 
Denote ${\cal L}(\tilde{\X}_n)=P_{\theta,n}$.
W.l.o.g., consider estimation of $\gamma_1=q\in\{1,\ldots,k\}$. Let $(G_n)_{n\in\N}$ be a consistent estimator of $\gamma_1$, meaning (because of discreteness) that 
\begin{equation} \label{eg1cons}
\lim_{n\to\infty}P_{\theta,n}\{G_n(\tilde{\X}_n)=q\}=1. 
\end{equation}
Consider now a parameter vector $\theta^*$ for which $\gamma_1=r\neq q$, but which is otherwise equal to $\theta$. Consistency implies
\begin{equation} \label{eg2cons}
\lim_{n\to\infty}P_{\theta^*,n}\{G_n(\tilde{\X}_n)=r\}=1. 
\end{equation}
As $\gamma_1$ determines the distribution of $\X_1$ but nothing else, for almost all $x\in\R:\ P_{\theta,n}(\tilde{\X}_n\mid\X_1=\x)=P_{\theta^*,n}(\tilde{\X}_n\mid\X_1=\x)$.
But then (\ref{eg1cons}) requires the existence of a set $B\subseteq \R:$
\[
\lim_{n\to\infty}P_{\theta,n}\{\X_1\in B\}=1,\ \forall \x\in B:\
\lim_{n\to\infty}P_{\theta,n}\{G_n(\tilde{\X}_n)=q\mid\X_1=\x\}=1,
\] 
and  (\ref{eg2cons}) requires
\[
\lim_{n\to\infty}P_{\theta^*,n}\{\X_1\in B\}=0.
\]
Such a set $B$ does not exist, because the distribution of $\X_1$ under both
$\theta$ and $\theta^*$ are fixed Gaussian distributions with nonzero density 
everywhere. Therefore $G_n$ cannot be consistent.
\subsection*{Proof of Theorem \ref{tidentpollard}}
W.l.o.g., consider estimation of $\gamma_1$. 
For $n\in\N$ let 
\[
G_n(\tilde{\X}_n)=\argmin_{j\in\{1,\ldots,k\}}\|\X_{1}-\m_{nj}\|^2
\]
be the $k$-means estimator of the cluster membership of $\X_1$ (in case of 
non-uniqueness of the $\argmin$, any tie breaking rule can be used). Because of 
(\ref{econslabels}) and the Theorem in \cite{Pollard81}, which applies because
of (\ref{epollardass}), 
\[
P^*\left\{\lim_{n\to\infty}T_n^m(\tilde{\X}_n)=(\bmu^*_1,\ldots,\bmu^*_k)\right\}=1.
\]
$\X_1$ is fixed once it is observed and does not change for $n\to\infty$. Unless $\exists j\neq l\in\{1,\ldots,k\}:\ \|\X_1-\bmu^*_j\|^2=\|\X_1-\bmu^*_l\|^2$ (which happens with probability 0 due to (\ref{econtass})), 
\[
\gamma_1=\argmin_{j\in\{1,\ldots,k\}}\|\X_1-\bmu^*_j\|^2, \mbox{ and } 
\epsilon=\min_{j\neq \gamma_1}(\|\X_1-\bmu^*_j\|^2-\|\X_1-\bmu^*_{\gamma_1}\|^2)>0.
\]
For $n$ large enough, with probability 1, therefore 
$\|\X_1-\m_{\gamma_1n}\|^2=\min_{j\in\{1,\ldots,k\}}\|\X_1-\m_{jn}\|^2$, and
$G_n(\tilde{\X}_n)=\gamma_1$. Therefore $G_n$ is consistent for $\gamma_1$,
and $\gamma_1$ as well as all $\gamma_i$ are empirically identifiable.

\bibliographystyle{chicago}

\end{document}